\documentclass[11pt]{article}

\usepackage{amsmath}
\usepackage{amsthm}
\usepackage{amssymb}
\usepackage{enumerate}
\usepackage{verbatim}
\usepackage{fullpage}

\newtheorem{theorem}{Theorem}[section]
\newtheorem{lemma}[theorem]{Lemma}
\newtheorem{prop}[theorem]{Proposition}
\newtheorem{corollary}[theorem]{Corollary}
\newtheorem*{conjecture}{Conjecture}
\newtheorem*{theorem1}{Theorem 1.1}

\theoremstyle{definition}
\newtheorem*{definition}{Definition}

\theoremstyle{remark}
\newtheorem*{remark}{Remark}

\newcommand{\vol}{\operatorname{vol}}
\newcommand{\Dist}{d}
\newcommand{\Grp}{\operatorname{Gr}_p}

\newcommand{\Zp}{\mathbb{Z}_p}
\newcommand{\Bset}{\{0\}\cup\{1/b \mid b\in \mathbb{Z}^+\}}

\begin{document}

\title{Nathanson heights in finite vector spaces}
\author{Joshua Batson \thanks{
\texttt{email: joshua.batson@yale.edu}}}
\date{}
\maketitle

\begin{abstract}
Let $p$ be a prime, and let $\mathbb{Z}_p$ denote the field of integers modulo $p$.  The \emph{Nathanson height} of a point $v \in \mathbb{Z}_p^n$ is the sum of the least nonnegative integer representatives of its coordinates.  The Nathanson height of a subspace $V \subseteq \mathbb{Z}_p^n$ is the least Nathanson height of any of its nonzero points.  In this paper, we resolve a conjecture of Nathanson [M. B. Nathanson, Heights on the finite projective line, International Journal of Number Theory, to appear], showing that on subspaces of $\mathbb{Z}_p^n$ of codimension one, the Nathanson height function can only take values about $p, p/2, p/3, \dots.$  We show this by proving a similar result for the coheight on subsets of $\mathbb{Z}_p$, where the \emph{coheight} of $A \subseteq \mathbb{Z}_p$ is the minimum number of times $A$ must be added to itself so that the sum contains 0.  We conjecture that the Nathanson height function has a similar constraint on its range regardless of the codimension, and produce some evidence that supports this conjecture.
\end{abstract}

\section{Introduction}

	Let $\mathbb{Z}_p$ denote the field of integers modulo $p$.  Let $(a_1,\dots,a_n) \in \mathbb{Z}^n$, where $n \geq 1$ and $0 \leq a_i < p$, and let $v = (\bar{a}_1,\dots,\bar{a}_n) \in \mathbb{Z}_p^n$, where $\bar{a} \in \mathbb{Z}_p$ denotes the residue class of $a \in \mathbb{Z}$.  The \emph{Nathanson height} of $v$ is $h_p(v) = a_1 + \dots + a_n$.  Let $V$ be a subspace of $\mathbb{Z}_p^n$.  The \emph{Nathanson height} of $V$ is 
	$$h_p(V) = \min \{ h_p(v) \mid v \in V\setminus \{0\} \}  =  \min \Big \{ \sum_{i=1}^m a_i \, \Big \vert \, 0 \leq a_i < p, (\bar{a}_1,\dots,\bar{a}_n) \in V\setminus \{0\} \Big \}.$$

	In \cite{N1}, Nathanson and Sullivan considered the properties of the Nathanson height function when $V$ is a one-dimensional subspace.  In this case, each line $V \subset \mathbb{Z}_p^n$ can be identified with a unique point $a$ in the projective space $\mathbb{P}^{n-1}$.  They proved $h_p(a) \leq n p /2$, then investigated the range of the Nathanson height function when $n = 2$.  In particular, they proved that if $h_p(a)$ is less than $p$, then it is in fact at most $(p+1)/2$.  In \cite{N2}, Nathanson extended this argument to show that if $h_p(a)$ is less than $(p+1)/2$, then it is at most $(p+4)/3$, that is, $h_p(a)$ is either $p$, about $p/2$, or at most roughly $p/3$.  He conjectured that for sufficiently large primes, this pattern continues: for fixed $b_0$ and sufficiently large $p$, if $a\in \mathbb{P}^1$, then either $0 \leq h_p(a)-p/b < c_b$ for some $b \leq b_0$ or $h_p(a) < p/b_0$, where each $c_b$ is a constant depending only on $b$.  In \cite{B1}, O'Bryant proved a weaker statement, that as $p$ tends to infinity, the set $\{\frac{1}{p}h_p(a) \mid a \in \mathbb{P}^1\}$ converges to $\Bset$.

	In this paper, we study the range of the Nathanson height function over $m$-dimensional subspaces of $\Zp^n$.  When $m = n-1$, that is, when the subspaces have codimension one, we obtain the following result:
\begin{theorem}\label{coheightvectorthm}
Let $p$ be a prime, let $n \geq 2$, and let $V$ be a subspace of $\mathbb{Z}_p^n$ of codimension one.  Then either
$$0 < h_p(V) < 3p^{(n-1)/n} + 1$$
or 
$$0 \leq h_p(V) - \frac{p}{b} < b$$
for some $1 \leq b \leq \left \lceil p^{1/n} \right \rceil.$
\end{theorem}
We also show that these inequalities can be satisfied, so the Hausdorff distance between $\{ \frac{1}{p} h_p(V) \mid V \subset \Zp^n \mbox{ has codimension one} \} $ and $\Bset$ is $O(p^{-1/n})$, where the Hausdorff distance between two sets $A$ and $B$ is $d(A,B) := \max \{\sup_{b \in B} \inf_{a \in A} |a-b|, \sup_{a \in A} \inf_{b \in B} |a-b| \}$.
	This result in the case $n=2$ proves Nathanson's conjecture.  We obtain the above results by transforming the problem of determining the Nathanson heights of subspaces of $\mathbb{Z}_p^n$ of codimension $m$ into a problem about sumsets in $\mathbb{Z}_p^m$: given $A \subseteq \Zp^m$, what is the least positive integer $k$ such that the $k$-fold sum $A+\dots+A$ contains 0?

	A brief outline of the paper follows.  In Section 2, we introduce notation and give a tight upper bound on the Nathanson height function.  In Section 3, we develop the link to sumsets, and in Section 4, we prove Theorem \ref{coheightvectorthm}.  In Section 5, we give a technique for calculating some values of the Nathanson height function for subspaces of codimension greater than one.  In Section 6, we show that that those subspaces with large Nathanson height in proportion to  $p$ are rare: for any $\varepsilon >0$ and for sufficiently large primes $p$, almost all $m$-dimensional subspaces of $\mathbb{Z}_p^n$ have Nathanson height less than $\varepsilon p$.  In Section 7, we discuss some open questions related to the Nathanson height function and give a conjecture on the range of the Nathanson height function for arbitrary codimension.

\section{Notation and upper bounds}

We begin our study of the range of the Nathanson height function by introducing some notation.  Let $p$ be a prime and let $1 \leq m < n$.  The Grassmannian $\operatorname{Gr}_p(n,m)$ is the set of $m$-dimensional subspaces of $\Zp^n$.  Let $H_p(n,m) = \{h_p(V) \mid V \in \operatorname{Gr}_p(n,m) \}$ and $H_p^c(n,m) = H_p(n,n-m)$; so $H_p^c(n,m)$ is the set of Nathanson heights of subspaces of $\Zp^n$ of codimension $m$.  We can now show how the range of the Nathanson height function changes with the dimension of the subspaces or the ambient space.

\begin{lemma}\label{htrelations}
Let $p$ be a prime and $1 \leq m < n$.  Then the following conditions hold:
\begin{enumerate}[(a)]
\item $H_p(n,m) \subseteq H_p(n,m-1) \subseteq \cdots \subseteq H_p(n,1)$;
\item $H_p^c(n,m) \supseteq H_p^c(n-1,m) \supseteq \cdots \supseteq H_p^c(m+1,m) = H_p(m+1,1)$;
\item $\max(H_p^c(n,m)) =  \max(H_p^c(n-1,m)) = \cdots = \max(H_p^c(m+1,m)).$
\end{enumerate}
\end{lemma}
\begin{proof}
\begin{enumerate}[(a)]
\item Let $V \in \operatorname{Gr}_p(n,m)$, and let $v \in V$ have minimum Nathanson height.  Then $v \in V^\prime \subset V$ for some $(m-1)$-dimensional subspace $V^\prime$ of $V$ and has the minimum Nathanson height over that subspace, so $h_p(V^\prime) = h_p(V)$.  Thus, $H_p(n,m) \subseteq H_p(n,m-1)$, and the statement follows because $m$ was arbitrary.

\item Let $W$ be a subspace of $\Zp^{n-1}$ of codimension $m$, and let $V = \phi^{-1}(W)$ where $\phi(a_1,\dots,a_n) = (a_1,\dots,a_{n-2},a_{n-1}+a_n).$  Then $V$ is a subspace of $\Zp^n$ of codimension $m$.  If $v \in V$ and $v_{n-1} + v_n < p$, then $h_p(v) = h_p(\phi(v))$, while if $v_{n-1} + v_n \geq p$, then $h_p(v)=h_p(\phi(v))+p$.  But any $w \in W$ has some preimage $v$ with $h_p(v) = h_p(w)$---just append $0$ to $w$ to get such a $v$.  So $h_p(V) = h_p(W)$.  Since $W$ was an arbitrary element of $\operatorname{Gr}_p(n-1,n-1-m)$, we have $H_p^c(n-1,m) \subseteq H_p^c(n,m)$, and the statement follows because $n$ was arbitrary.

\item By (b), it suffices to show that $\max(H_p^c(n,m)) \leq \max(H_p^c(m+1,m))$.  Let $V$ be a subspace of $\Zp^n$ of codimension $m$.  Then there is some nonzero $v \in V$ with $0$ in its first $n-m-1$ coordinates.  Let $v^\prime \in \mathbb{Z}_p^{m+1}$ have coordinates equal to the $m+1$ remaining coordinates of $v$, and let $V^\prime$ be the one-dimensional subspace of $\mathbb{Z}_p^{m+1}$ spanned by $v^\prime$.  Note that $V^{\prime}$ consists of multiples of $v^\prime$, every multiple of $v$ is in $V$, and $h_p(kv) = h_p(kv^\prime)$ for any $k \in \mathbb{Z}$.  Therefore, $h^c_p(V) \leq h^c_p(V^\prime)$.  Since $V$ was chosen arbitrarily and $V^\prime$ has codimension $m$ in $\Zp^{m+1}$, we have $\max(H_p^c(n,m)) \leq \max(H_p^c(m+1,m))$. \qedhere
\end{enumerate} \end{proof}

As the above lemma suggests, the codimension of a subspace of $\Zp^n$ greatly constrains its Nathanson height.  For example, the Nathanson height of a subspace can be bounded above by a function depending only on its codimension.

\begin{lemma}\label{coheightupperbound}
Let $p$ be a prime, and let $1 \leq m < n$.  Then 

$$\max(H_p^c(n,m)) = 
\begin{cases}
\frac{(m+1)p}{2}, \mbox{ if $m$ is odd};\\
\frac{mp}{2}+1, \mbox{ if $m$ is even}.
\end{cases}
$$
\end{lemma}

\begin{proof}
By part (c) of Lemma \ref{htrelations}, we need only consider the case $n = m+1$.  Then $H_p^c(n,m) = H_p(n,1)$.  Let $V$ be a one-dimensional subspace of $\mathbb{Z}_p^n$ generated by $v \in \Zp^n$.  We have $h_p(-v) + h_p(v) \leq np$, since each nonzero coordinate of $v$ contributes $p$ to the sum, so either $h_p(v)$ or $h_p(-v)$ is at most $np/2$.  Then $h_p(V) \leq np/2$, so $\max(H_p(n,1)) \leq np/2 = (m+1)p/2$.  If $m$ is odd, then $n = m+1$ is even and this bound is tight:  the subspace generated by the point $(1,p-1,1,p-1,\dots,1,p-1) \in \mathbb{Z}_p^{m+1}$ has Nathanson height $(m+1)p/2$.  

Now suppose $m$ is even and $h_p(V) > mp/2+1 = (n-1)p/2+1$.  Let $v \in V$.  Then $h_p(v) \geq (n-1)p/2+2$ and $$np - h_p(v) \geq h_p(-v) \geq \frac{(n-1)p}{2}+2,$$ so
	$$\frac{(n-1)p}{2}+2 \leq h_p(v) \leq \frac{(n-1)p}{2} + p-2.$$
There are $p-3$ values in this range, but $V$ contains $p-1$ nonzero points.  So there exist distinct $v_1,v_2 \in V$ with $h_p(v_1) = h_p(v_2)$.  Let $v^\prime =  v_1-v_2$, so 
$$h_p(v^\prime) = h_p(v_1-v_2) \equiv h_p(v_1) - h_p(v_2) \equiv 0 \pmod{p}.$$
Since $n=m+1$ is odd and either $h_p(v^\prime) \leq np/2$ or $h_p(-v^\prime) \leq np/2$, we have 
	$$h_p(V) \leq \min(h_p(v^\prime),h_p(-v^\prime)) \leq \frac{(n-1)p}{2} =  \frac{mp}{2},$$
which is a contradiction.  Thus $h_p(V) \leq mp/2+1$ and $\max(H_p(n,1))  \leq mp/2+1$.  This bound is tight: the subspace generated by $(1,p-1,1,p-1,\dots,1,p-1,1) \in \mathbb{Z}_p^{m+1}$ has Nathanson height $mp/2+1$.
\end{proof}

\section{Coheight and sumsets}

We can transform the problem of the Nathanson heights of subspaces of $\mathbb{Z}_p^n$ of codimension $m$ into a problem about sumsets in $\mathbb{Z}_p^m$ as follows:  let $V$ be a subspace of $\mathbb{Z}_p^n$ of codimension $m$ with $1 \leq m < n$.  Then $V$ is the kernel of some $m \times n$ matrix $M$ with columns $a_j$.  So $b \in V$ when $Mb = 0$, that is, when $\sum \bar{b}_j a_j = 0$.  Therefore, there exists a point $b = (\bar{b}_1, \dots, \bar{b}_n) \in V$ with $h_p(b) = \sum_{i=1}^n b_j = k$ if and only if there exist exactly $k$ (not necessarily distinct) elements of $A = \{a_1,\dots, a_n\} \subset \mathbb{Z}_p^m$ that sum to zero with every element of $A$ occurring fewer than $p$ times in the sum.

	To formalize this duality between sumsets in $\Zp^m$ and the height, we require the following definitions.  Let $G$ be an abelian group.  For $A,B \subseteq G$, let $A+B$ denote $\{a+b \mid a \in A, b \in B \}$.  For $k$ a positive integer, let $kA$ denote the $k$-fold sum $A+\dots+A$.  We can now define a coheight on an arbitrary finite abelian group.
	
\begin{definition}
Let $G$ be a finite abelian group, and let $A \subseteq G$, $A \neq \emptyset$.  The \emph{coheight} of $A, h_G^c(A),$ is the least positive integer $k$ such that $0 \in kA$.
\end{definition}

In the case where $G = \mathbb{Z}_p^m$, the coheight provides a lower bound for the Nathanson height on subspaces of $\Zp^n$ of codimension $m$.  In what follows, we write $h_p(A)$ for $h_{\mathbb{Z}_p^m} (A)$.

\begin{prop}\label{dualityprop}
Let $V$ be a subspace of $\Zp^n$ of codimension $m$ with $1 \leq m < n,$ so $V$ is the kernel of some $m \times n$ matrix $M$ with columns $a_j$.  Let $A = \{a_1,\dots,a_n\} \subseteq \Zp^m$.  Then $h_p^c(A) \leq h_p(V)$.  Moreover, if $h_p^c(A) < p$, then $h_p^c(A) = h_p(V)$.
\end{prop}
\begin{proof}
Let $b = (\bar{b}_1, \dots, \bar{b}_n) \in V$ have minimum Nathanson height, $h_p(b) = \sum b_j = h_p(V)$.  Then $\sum b_j a_j = 0 \in h_p(V)A$, so $h_p^c(A) \leq h_p(V)$.

If $h_p^c(A) < p$, then there exist nonnegative integers $b_j$ with $\sum b_j = h_p^c(A) <p$ and $\sum b_j a_j = 0$.  Thus $0 \leq b_j <p$ and $b = (\bar{b}_1,\dots,\bar{b}_n) \in V$, so 
$$h_p(V) \leq h_p(b) = \sum b_j = h_p^c(A).$$
Therefore, $h_p(V)=h_p^c(A)$.
\end{proof}

Equality does not hold in general---we have $h_p(A) \leq p$ since $0 \in pA$, while $h_p(V)$ can be as large as $mp/2$.  In the case where $V$ has codimension one, though, Lemma \ref{coheightupperbound} gives $h_p(V) \leq p$, so $h_p(V) = h_p^c(A)$.  Thus we can use the coheight of subsets of $\mathbb{Z}_p$ to calculate the Nathanson height of subspaces of $\mathbb{Z}_p^n$ of codimension one.  Before narrowing our focus to $\mathbb{Z}_p$, we define a variation on the coheight function.

\begin{definition}
Let $G$ be a finite abelian group, and let $A$ be a nonempty subset of $G$ not contained in any proper subgroup.  The \emph{width} of $A, w_G(A),$ is the least positive integer $k$ such that every element of $G$ can be expressed as the sum of at most $k$ elements of $A$, that is, $\bigcup_{1 \leq k^\prime \leq k} k^\prime A = G$.
\end{definition}

The width is clearly an upper bound on the coheight.  We can in turn give an upper bound on the width $w_G(A)$ in terms of the sizes of $A$ and $G$.  We require the following result by Kneser in our proof.  (A straightforward proof is given in \cite{T1}, p. 200.)

\begin{theorem}[Kneser \cite{K1}] \label{kneser}
Let $G$ be a finite abelian group, and let $A,B \subseteq G$.  Then 
$$|A + B| \geq|A|+|B|-|\operatorname{Sym}_G(A+B)|,$$
where $\operatorname{Sym}_G(X)$ is the subgroup $\{h \in G \mid X+h = X\}$.
\end{theorem}

\begin{lemma}\label{widthlemma}
Let $G$ be a finite abelian group, and let $A$ be a nonempty subset of $G$ not contained in any proper subgroup.  If $0 \in A$, then $w_G(A) \leq \left \lceil \frac{2|G|}{|A|} \right \rceil - 1$, while if $0 \notin A$, then $w_G(A) \leq \left \lceil \frac{2|G|}{|A|+1} \right \rceil$.
\end{lemma}
\begin{proof}

We first assume $0 \in A$.  In that case, $x \in G$ can be expressed as the sum of at most $k$ elements of $A$ if and only if it can be expressed as the sum of exactly $k$ elements of $A$, that is, $x \in kA$.  Then $w_G(A)$ is the least positive integer $k$ such that $kA = G$.

	We proceed by induction on the number of proper subgroups of $G$.  Suppose $G$ has no proper nontrivial subgroups, and let $B \subseteq G$.  Then Theorem \ref{kneser} gives $|A + B| \geq |A|+|B|-1$ unless $A+B = G$.  Setting $B = kA$ and inducting gives $|nA| \geq n|A|-n+1$ unless $nA = G$.  If $n \geq  \frac{|G|-1}{|A|-1} $, the previous inequality implies $nA = G$.  Note $\frac{2|G|}{|A|} - 1\geq \frac{|G|-1}{|A|-1}$, since the inequality holds when $|G| = |A|$ and is strengthened by $\frac{2}{|A|} - \frac{1}{|A|-1} \geq 0$ if we increment $|G|$.  Thus, $w_G(A) \leq \left \lceil \frac{2|G|}{|A|} \right \rceil - 1$.
	
	Now consider an arbitrary finite abelian group $G$.  Since $\operatorname{Sym}_G((n-1)A) \subseteq \operatorname{Sym}_G(nA)$, we can inductively show that $|nA| \geq n|A| - (n-1) |\operatorname{Sym}_G(nA)|$.  Let $n = \left \lceil \frac{2|G|}{|A|} \right \rceil - 1$.  Then $H = \operatorname{Sym}_G(nA)$ cannot be trivial, since $|G| \geq |nA|$.  If $H = G$, then $nA = G$ and $w_G(A) \leq \left \lceil \frac{2|G|}{|A|} \right \rceil - 1$ as desired.  So suppose $H$ is a nontrivial proper subgroup of $G$.  Let $A^\prime = (A+H)/H$, so $(nA)/H = nA^\prime$.  Since $nA$ is a union of cosets of $H$, we have that $nA^\prime = G/H$ implies $nA = G$.  Thus $w_G(A) \leq \max(w_{G/H}(A^\prime),n)$.  Since $0 \in A^\prime$, the inductive hypothesis gives
	$$w_{G/H}(A^\prime) \leq \left \lceil \frac{2|G/H|}{|A^\prime|} \right \rceil - 1 \leq \left \lceil \frac{2|G|/|H|}{|A|/|H|} \right \rceil - 1 = \left \lceil \frac{2|G|}{|A|} \right \rceil - 1 = n.$$
Thus $w_G(A) \leq \max(n,n) =  \left \lceil \frac{2|G|}{|A|} \right \rceil- 1$, and the induction is complete.

We now consider the case $0 \notin A$.  Any $x \in G \setminus \{0\}$ can be expressed as the sum of at most $w_G(A \cup \{0\})$ elements of $A \cup \{0\}$, and hence of at most $w_G(A \cup \{0\})$ elements of $A$.  Since $0 = (-a) + a$ for $a \in A$, we can express $0$ as the sum of at most $w_G(A \cup \{0\})+1$ elements of $A$.  Therefore, 
$$w_G(A) \leq w_G(A \cup \{0\}) + 1 \leq \left \lceil \frac{2|G|}{|A \cup \{0\}|} \right \rceil - 1 + 1 \leq \left \lceil \frac{2|G|}{|A|+1} \right \rceil,$$
which completes the proof.
\end{proof}

Without more information about the structure of the set $A$ or the underlying group $G$, it is difficult to say more about the values of the coheight or the width.  One thing we can say, though, is that these values are invariant under automorphisms of $G$.

\begin{lemma}\label{invariance lemma}
Let $G$ be a finite abelian group, and let $A = \{a_1,\dots,a_n\}$ be a nonempty subset of $G$ not contained in any proper subgroup.  Let $\sigma$ be an automorphism of G.  Then
$$h^c_G(A) = h^c_G(\sigma(A)) \mbox {  and  }  w_G(A) = w_G(\sigma(A)).$$
In addition, if $0 \notin A$, then $w_G(A) = w_G(\{-a_1, a_2 - a_1,\dots,a_n-a_1\})$.

\end{lemma}
\begin{proof}
The first statement follows because automorphisms preserve addition and send both $0$ and $G$ to themselves.

Now let $g$ be an arbitrary element of $G$, and let $k = w_G(A)$.  We can write $g = k_1a_1 + \dots + k_n a_n$, with $k_1 + \dots + k_n = k^\prime \leq k$.  So 
\begin{align*}
g - ka_1 &= (k - k^\prime)(-a_1) + k_1 (a_1 - a_1) + k_2(a_2-a_1)+ \dots + k_n (a_n - a_1) \\
	      &= k_1^\prime(-a_1) + k_2^\prime (a_2 - a_1) + \dots + k_n^\prime (a_n - a_1),
\end{align*}
where $k_1^\prime = k - k^\prime$ and $k_i^\prime = k_i$ for $i \geq 2$.  Let $A^\prime = \{-a_1, a_2 - a_1,\dots,a_n-a_1\}$, so $g$ can be written as the sum of $(k - k^\prime) + k_2 + \dots + k_n = k - k_1 \leq k$ elements of $A^\prime$.  Since $G - ka_1 = G$ and $g$ was chosen arbitrarily, $w_G(A) \leq w_G(A^\prime).$ Since $0 \notin A$, the above procedure is reversible, so $w_G(A^\prime) \leq w_G(A)$.  Therefore, $w_G(A) = w_G(A^\prime) = w_G(\{-a_1, a_2 - a_1,\dots,a_n-a_1\})$, as desired.
\end{proof}

If $A \subset G = \mathbb{Z}_p$ and $k \in \mathbb{Z}_p\setminus \{0\}$, the above lemma gives
$h_G^c(A) = h_G^c(k \cdot A)$ and $w_G(A) = w_G(k \cdot A)$, where $k \cdot A = \{ ka \mid a \in A \}$.  We use this fact in the next section.

%%%%%
\section{Codimension one and sumsets in $\mathbb{Z}_p$}

In this section, we calculate the range of the coheight and width functions on subsets of $\mathbb{Z}_p$.  As a corollary, this gives us the range of the Nathanson height function on subspaces of $\mathbb{Z}_p^n$ with codimension 1.  In particular, we prove the conjecture of Nathanson on the range of the Nathanson height function in $\mathbb{Z}_p^2.$

While it is difficult to calculate the coheight and width in general, we can order the elements of $\mathbb{Z}_p$ by their representatives in $\{0,1,\dots,p-1\}$.  This allows us to give the following bounds on the coheight and the width functions on $\mathbb{Z}_p$.  In this section, we write $h_p^c$ for $h_{\Zp}^c$ and $w_p$ for $w_{\Zp}$.

\begin{lemma}\label{htboundlemma}
Let $0<a_1 < \cdots < a_n < p$ with $\gcd(a_1,\dots,a_n) = 1$, and let $A = \{\bar{a}_1,\dots, \bar{a}_n\} \subset \mathbb{Z}_p$.  Then
$$\frac{p}{a_n} \leq h^c_p(A) \leq w_p(A) < \frac{p}{a_n} + a_n+a_{n-1}.$$
Moreover, if $a_{n-1}(a_n-1) < p$, then we have the stronger statement $h_p^c(A) < \frac{p}{a_n} + a_n.$
\end{lemma}
\begin{proof}
	The first two inequalities follow from elementary considerations:  we have $\frac{p}{a_n} \leq h_p^c(A)$ because if $0 < k < \frac{p}{a_n}$, then $kA \subset \{k\bar{a}_1,k\bar{a}_1+1,\dots,k\bar{a}_n\}$, which does not contain $0$.  The second inequality, $h_p^c(A) \leq w_p(A)$, is trivial.  If $|A| = 1$, then $a_1 = 1$, so $h_p^c(A) = w_p(A) = p$ and the remaining inequalities hold.

	Suppose $|A|>1$.  Let $A^\prime \subseteq \mathbb{Z}_{a_n}$ be the set of residues of $0,a_1,\dots,a_{n-1}$ in $\mathbb{Z}_{a_n}$.  Then $A^\prime$ is not contained in a proper subgroup of $\mathbb{Z}_{a_n}$, since $\gcd(a_1,\dots,a_n) = 1$.  By Lemma \ref{widthlemma},
	$$w_{\mathbb{Z}_{a_n}}(A^\prime) \leq \left \lceil \frac{2a_n}{|A^\prime|} \right \rceil - 1\leq a_n - 1,$$
so any element of $\mathbb{Z}_{a_n}$ can be written as the sum of at most $a_n-1$ elements of $A^\prime$.
	Let $\bar{x} \in \mathbb{Z}_p$, and let $x^\prime$ be the unique integer congruent to $x$ in $\{a_{n-1} (a_n-1),a_{n-1} (a_n-1)+1,\dots,a_{n-1}(a_n-1) + p- 1\}$.  There exists a set of at most $a_n-1$ (not necessarily distinct) elements of $\{0,a_1,\dots,a_{n-1}\}$ whose sum $s$ is congruent to $x^\prime$ modulo $a_n$.  Clearly $s$ lies between $0$ and $a_{n-1}(a_n-1)$, so $0 \leq x^\prime - s \leq p+a_{n-1}(a_n-1) - 1$.  Since $x^\prime = \frac{x^\prime - s}{a_n}a_n + s$ and $x^\prime$ has residue $\bar{x}$ in $\mathbb{Z}_p$, we have expressed $\bar{x}$ as the sum of at most 
	$$\frac{p+a_{n-1}(a_n-1) - 1}{a_n} + a_n < \frac{p}{a_n} + a_{n-1}+ a_n$$
	 elements of $A$.  But $\bar{x}$ was an arbitrary element of $\mathbb{Z}_p$, so $w_p(A) < \frac{p}{a_n} + a_{n-1}+ a_n$.
	
	Suppose that $a_{n-1}(a_n-1) < p$.  Setting $\bar{x} = \bar{0}$ in the above argument, we get that the unique $x^\prime \in \{a_{n-1} (a_n-1),\dots,a_{n-1} (a_n-1)+p - 1\}$ congruent to $\bar{0}$ is precisely $p$.  This gives $0 \leq x^\prime - s < p$, which in turn yields a set of fewer than $\frac{p}{a_n}+ a_n$ elements of $A$ that sum to $\bar{0}$.  Then $h^c_p(A) < \frac{p}{a_n} + a_n.$
\end{proof}

The above lemma gives good bounds on the coheight and width when the elements of $A$ are small.  The following lemma tells us how small we can make them using the automorphisms of $\mathbb{Z}_p$.  For notational convenience, we define  $|\bar{x}|_p = \min(x, p - x)$ for any $0 \leq x < p$.

\begin{lemma} \label{pigeonlemma}  Let $p$ be prime, let $n>1$, and let $A  = \{\bar{a}_1,\dots, \bar{a}_n\} \subset \mathbb{Z}_p$.  Then for some $k \in \{1,\dots,p-1\}$, we have $|k \bar{a}_1|_p,\dots,|k \bar{a}_n|_p \leq \left \lceil p^{ (n-1)/n} \right \rceil$.
\end{lemma}

\begin{proof}  Let $s =\left \lceil p^{ (n-1)/n} \right \rceil $, and let $a = (\bar{a}_1,\dots,\bar{a}_n) \in \mathbb{Z}_p^n$.  Using a pigeonhole argument, we will show that for some $k \in \{1,\dots,p-1\}$, $k a$ lies in $\{-\bar{s}, -\bar{s}+1, \dots, \bar{s}\}^n \subset \mathbb{Z}_p^n$.  Let $q = \left \lfloor \frac{p}{s} \right \rfloor$.  We can then partition $\mathbb{Z}_p$ into $q$ intervals of length $s$ or $s+1$.  This division of $\mathbb{Z}_p$ into $q$ intervals yields a division of $\mathbb{Z}_p^n$ into $q^{n}$ $n$-dimensional boxes with sides $s$ or $s+1$.  We have 

$$q^{n} = \left \lfloor \frac{p}{s} \right \rfloor^{n} < \left (\frac{p}{s} \right )^{n} < \left (\frac{p}{p^{ (n-1)/n}}\right )^{n} = p.$$

	Now consider the $p$ points $0,a, 2a, \dots, (p-1)a \in \mathbb{Z}_p^{n}$, which lie in the $q^n \leq p-1$ boxes.  By the pigeonhole principle, some two points, say $k_0 a = (k_0\bar{a}_1,\dots,k_0\bar{a}_n)$ and $k_1 a =(k_1\bar{a}_1,\dots,k_1\bar{a}_n)$, lie in the same box.  Since each side of the box is at most $s+1$, $|k_1\bar{a}_i - k_0 \bar{a}_i|_p \leq s$ for $1 \leq i \leq n$.  Setting $k = |k_1 - k_0|$ gives the desired result.
\end{proof}

By combining the above lemmas, we can bound the range of the coheight function on subsets of $\mathbb{Z}_p$ with small cardinality.

\begin{theorem}\label{coheightthm}
Let $p$ be a prime, and let $A \subseteq \mathbb{Z}_p$ be a subset of size $n \geq 1$.  Then either 
\begin{equation}\label{eq1}
0 < h^c_p(A) < 3p^{(n-1)/n} + 1
\end{equation}
or 
\begin{equation}\label{eq2}
0 \leq h^c_p(A) - \frac{p}{b} < b
\end{equation}
for some $n \leq b \leq \left \lceil p^{1/n} \right \rceil.$
\end{theorem}
\begin{proof}
Let $A = \{\bar{a}_1,\dots,\bar{a}_n\}$.  If $n=1$, then either $\bar{a}_1 = 0$, in which case $h^c_p(A) = 1$ and \eqref{eq1} holds, or $\bar{a}_1 \neq 0$, in which case $h^c_p(A) = p$ and \eqref{eq2} holds for $b = 1$.  So we assume $n \geq 2$.

By Lemma \ref{pigeonlemma}, there is some $k \in \{1,\dots,p-1\}$ such that $|k \bar{a}_1|_p,\dots,|k \bar{a}_n|_p \leq \left \lceil p^{ (n-1)/n} \right \rceil$.  Since $h_p^c(A)$ does not change if we permute the $\bar{a}_i$ or multiply them by some such $k$, we can assume without loss of generality that $0 \leq |\bar{a}_1|_p \leq \dots \leq |\bar{a}_n|_p \leq \left \lceil p^{ (n-1)/n} \right \rceil$.  We can also assume $\gcd(|\bar{a}_1|_p,\dots,|\bar{a}_n|_p) = 1$, since if not, we multiply $A$ by the inverse of the gcd modulo $p$.  Let $a_i$ be the least nonnegative integer representative of $\bar{a}_i$ and assume $a_1 = |\bar{a}_1|_p$, since if not, we multiply $A$ by $-1$.  If $a_1 = 0$,  then $h_p^c(A) = 1$ and \eqref{eq1} holds.  If $p -  \left \lceil p^{ (n-1)/n} \right \rceil \leq a_i < p $ for some $2 \leq i \leq n$, then $a_1 \bar{a}_i + (p-a_i) \bar{a}_1 = 0$ and $0 \in (a_1 + p - a_i)A$.  In that case, 
$$h_p^c(A) \leq a_1 + p - a_i \leq |\bar{a}_1|_p + |\bar{a}_i|_p \leq 2p^{(n-1)/n} + 2 < 3p^{(n-1)/n} + 1$$
and \eqref{eq1} holds.

	If $0 < a_i \leq \left \lceil p^{ (n-1)/n} \right \rceil$ for all $2 \leq i \leq n$, then we have $0 < a_1 < \cdots < a_n$ and $\gcd(a_1,\dots,a_n) = 1$.  If $1 \leq a_n \leq \left \lceil p^{1/n} \right \rceil$, then $a_{n-1}(a_n-1) < p$.  By Lemma \ref{htboundlemma},
	$$\frac{p}{a_n} \leq  h_p^c(A) < \frac{p}{a_n} + a_n.$$
Let $b = a_n$ and subtract $1/b$ to give \eqref{eq2}.

In the final case, we have $\left \lceil p^{1/n} \right \rceil < a_n \leq \left \lceil p^{(n-1)/n} \right \rceil.$  Then Lemma \ref{htboundlemma} yields
	$$0 < h_p^c(A) \leq \frac{p}{a_n} + a_n+a_{n-1} < \frac{p}{p^{1/n}} + 2\left \lceil p^{(n-1)/n} \right \rceil - 1 < 3p^{(n-1)/n} + 1,$$
as desired.
\end{proof}

We now use the duality between Nathanson heights of subspaces of $\mathbb{Z}_p^n$ of codimension one and coheights of subsets of $\mathbb{Z}_p$.  By Proposition \ref{dualityprop}, if $V$ is the subspace of $\mathbb{Z}_p^n$ of codimension one that is the kernel of the $1\times n$ matrix $(\bar{a}_1 \dots \bar{a}_n)$ and $A = \{\bar{a}_1,\dots,\bar{a}_n\} \subset \mathbb{Z}_p$, then $h_p(V) = h_p^c(A)$.  The size of $A$ is equal to the number of distinct $\bar{a}_i$, which can be any number from $1$ to $n$.  Thus $H_p^c(n,1) = \{h_p(A) \mid A \subset \mathbb{Z}_p, 1 \leq |A| \leq n \}$.  Then Theorem \ref{coheightvectorthm} follows from Theorem \ref{coheightthm}.

\begin{theorem1}  
Let $p$ be a prime, let $n \geq 2$, and let $V$ be a subspace of $\mathbb{Z}_p^n$ of codimension one.  Then either
$$0 < h_p(V) < 3p^{(n-1)/n} + 1$$
or 
$$0 \leq h_p(V) - \frac{p}{b} < b$$
for some $1 \leq b \leq \left \lceil p^{1/n} \right \rceil.$
\end{theorem1}

When $n = 2$, the statement can be improved slightly by eliminating the case $\left \lceil p^{1/n} \right \rceil < a_n \leq \left \lceil p^{(n-1)/n} \right \rceil$, because $1/n = (n-1)/n = 1/2$, and the case $|\bar{a}_1|_p = |\bar{a}_i|_p = \left \lceil \sqrt{p} \right \rceil$, because their $\gcd$ is greater than 1.  This consideration yields the following result, which proves Nathanson's conjecture in \cite{N2} on the range of the Nathanson height function on the projective line.

\begin{corollary}\label{nathansoncor}
Let p be a prime and let $V$ be a line in $\mathbb{Z}_p^2$.  Then either
$$0 < h_p(V) < 2\sqrt{p} + 1$$
or 
$$0 \leq h_p(V) - \frac{p}{b} < b$$
for some $1 \leq b \leq \left \lceil \sqrt{p} \right \rceil.$
\end{corollary}

Theorem \ref{coheightvectorthm} also implies that as $p \rightarrow \infty$, the sets $\frac{1}{p} H_p^c(n,1) = \{\frac{h_p(V)}{p} \mid V \in \operatorname{Gr}_p(n,n-1) \}$ converge to $\Bset$ with respect to the Hausdorff distance.

\begin{theorem}\label{asymheight}
Let $p$ be a prime and $n \geq 2$.  The Hausdorff distance between $\frac{1}{p} H_p^c(n,1)$ and $\Bset$ is $O(p^{-1/n})$.
\end{theorem}
\begin{proof}
Let $B =\Bset $ and let $x \in \frac{1}{p} H_p^c(n,1)$.  By Theorem \ref{coheightvectorthm}, $x$ is either at most $3p^{-1/n} + p^{-1}$ away from $0$ or at most $\left \lceil p^{1/n} \right \rceil p^{-1} \leq 2p^{-1/n}$ away from $1/b \in B$.  Thus every element of  $\frac{1}{p} H_p^c(n,1)$ is at most $O(p^{-1/n})$ from some element of $B$.  Conversely, let $x \in B$.  If $x \leq \left \lceil p^{-1/n} \right \rceil$, then $x$ is at most $\left \lceil p^{-1/n} \right \rceil$ away from $\frac{1}{p}h^c_p(\{1,p-1\}) = \frac{2}{p} \in \frac{1}{p} H_p^c(n,1)$.  Otherwise, $1 \geq x = \frac{1}{b} > \left \lceil p^{-1/n} \right \rceil$, so $1 \leq b < p^{1/n}$.  Then $x$ is within $p^{-(n-1)/n} \leq p^{-1/n}$ of $\frac{1}{p}h^c_p(\{1,b\}) \in \frac{1}{p} H_p^c(n,1)$ by Lemma \ref{htboundlemma}.  Thus every element of $B$ is at most $O(p^{-1/n})$ from some element of $\frac{1}{p} H_p^c(n,1)$.
\end{proof}

We can also bound the range of the width function on subsets of $\mathbb{Z}_p$.

\begin{theorem}\label{widththm}
Let $p$ be a prime, and let $A \subseteq \mathbb{Z}_p$ with $n = |A \setminus \{0\}| \geq 1$.  Then either 
\begin{equation}\label{eq3}
0 < w_p(A) < 5p^{(n-1)/n} + 2
\end{equation}
or 
\begin{equation}\label{eq4}
-1 \leq w_p(A) - \frac{p}{b} < 2b
\end{equation}
for some $n \leq b \leq \left \lceil p^{1/n} \right \rceil.$
\end{theorem}
\begin{proof}
Let $A = \{\bar{a}_1,\dots,\bar{a}_n\}$.  If $n=1$, then $w_p(A) = p$ and \eqref{eq2} holds with $b = 1$, so we assume $n \geq 2$.
	We first consider the case $0 \notin A$.  By Lemma \ref{pigeonlemma}, there is some $k \in \{1,\dots,p-1\}$ such that $|k \bar{a}_1|_p,\dots,|k \bar{a}_n|_p \leq \left \lceil p^{ (n-1)/n} \right \rceil$.  Let $k \bar{a}_i$ be the least of $k \bar{a}_1, \dots k \bar{a}_n$, where the elements are ordered by their representatives in $[-(p-1)/2,(p-1)/2]$.  If that representative of $k \bar{a}_i$ is nonnegative, then let $\bar{a}_j^\prime = k\bar{a}_j$ for every $j$.  Otherwise, let $\bar{a}_j^\prime = k \bar{a}_j - k \bar{a}_i$ for $j \neq i$, and let $\bar{a}_i^\prime = -k \bar{a}_i$.  By Lemma \ref{invariance lemma},  $w_p(A) = w_p(A^\prime)$, where $A^\prime = \{\bar{a}_1^\prime, \dots, \bar{a}_n^\prime\}$.  Relabel the elements of $A^\prime$ to get $0 < a_1^\prime < a_2^\prime < \dots < a_n^\prime \leq 2 \left \lceil p^{ (n-1)/n} \right \rceil$, where $a_j^\prime$ is the least nonnegative integer representative of $\bar{a}_j^\prime$.  As in the proof of Theorem \ref{coheightthm}, we can assume $\gcd(a_1^\prime,\dots,a_n^\prime) = 1$.  By Lemma \ref{htboundlemma},
$$ \frac{p}{a_n^\prime} \leq w_p(A) < \frac{p}{a_n^\prime} + a_n^\prime+a_{n-1}^\prime. $$

As in the proof of Theorem \ref{coheightthm}, we analyze two cases.  If $0 < a_n^\prime \leq \left \lceil p^{1/n} \right \rceil$, then \eqref{eq4} holds for $b = a_n^\prime$, while if $\left \lceil p^{1/n} \right \rceil < a_n^\prime \leq 2 \left \lceil p^{ (n-1)/n} \right \rceil$, then \eqref{eq3} holds.

If $0 \in A$, then we have $w_G(A)$ equals $w_G(A \setminus \{0\})-1$ or $w_G(A \setminus \{0\})$, depending on whether or not every nonzero element of $G$ can be obtained using at most $w_G(A \setminus \{0\})-1$ elements of $A$.  The statement then follows from the application of the first case to $A \setminus \{0\}$.
\end{proof}

We can now show that the sets $\{\frac{1}{p}h^c_p(A) \mid A \subset \mathbb{Z}_p\}$ and $\{\frac{1}{p}w_p(A) \mid A \subset \mathbb{Z}_p\}$ converge to the set $\Bset $ as $p \rightarrow \infty$.  Unlike in Theorem \ref{asymheight}, we allow $|A|$ to be arbitrary, even as large as $p$.  (For notational convenience, we set $h_p(\emptyset) = w_p(\emptyset) = w_p(\{0\}) = p$.)  We can bound the rate of convergence as follows.

\begin{corollary}
Let $p \geq 3$ be a prime.  Then the Hausdorff distance between $\Bset$ and each of $\{\frac{h^c_p(A)}{p} \mid A \subseteq \mathbb{Z}_p\}$ and $\{\frac{w_p(A)}{p} \mid A \subseteq \mathbb{Z}_p\}$ is $O(\frac{\log \log p}{\log p})$.
\end{corollary}
\begin {proof}
If $n < \frac{\log p}{\log \log p}$, then $p^{-1/n} < \frac{1}{\log p}$.  Since
$H_p^c(n,1) = \{h^c_p(A) \mid A \subset \mathbb{Z}_p, 1 \leq |A| \leq n\}$ and $d(\frac{1}{p}H_p^c(n,1), B) = O(p^{-1/n}),$ we have 
$$d\left (\left \{\frac{h^c_p(A)}{p}  \mid A \subset \mathbb{Z}_p, 1 \leq |A| \leq n\right \},B\right ) \leq O\left (p^{-1/n}\right ) \leq O \left (\frac{1}{\log p}\right ) \leq  O\left (\frac{\log \log p}{\log p}\right ).$$
But if $|A| \geq \frac{\log p}{\log \log p} > 0,$ then Lemma \ref{widthlemma} gives 
$$\frac{h_p^c(A)}{p} \leq \frac{w_p(A)}{p} \leq \frac{1}{p}\left \lceil \frac{2p}{|A|} \right \rceil < \frac{2 \log \log p}{\log p} + \frac{1}{p},$$
which is $O(\frac{\log \log p}{\log p})$ away from $0$.  Thus $d(\{\frac{1}{p}h^c_p(A) \mid A \subseteq \mathbb{Z}_p\},B)$ is $O(\frac{\log \log p}{\log p})$.
	The corresponding statement for the width follows from a similar argument using Theorem \ref{widththm}.
\end{proof}

%%%%%
\section{Codimension two, and beyond}
	In Section 4, we proved that the sets $\frac{1}{p}H_p^c(n,1)$ converge to $\Bset$.  For $m>2$, it is an open question whether or not the sets $\frac{1}{p}H_p^c(n,m)$ converge at all, much less to a particular set.  However, we can say something about the topological limit superior of $\frac{1}{p}H_p^c(n,m)$, where $x \in \limsup_{q \rightarrow \infty} \frac{1}{q} H^c_q(n,m)$ if there exist infinite sequences of primes $\{q_i\}$ and elements $x_i \in \frac{1}{q_i} H^c_{q_i}(n,m)$ such that the $x_i$ converge to $x$.  In this section, we give an infinite family of points contained in $\limsup_{q \rightarrow \infty} \frac{1}{q} H^c_q(n,m)$ for any $1 \leq m <n$.  In the case $m=2$, we show that
	$$\bigcup_p \frac{1}{p} H_p(2,1) \subset \limsup_{q \rightarrow \infty} \frac{1}{q} H_q^c(n,2). $$
In other words, the values of the Nathanson height function for lines in $\Zp^2$ generate limiting values of the Nathanson height function for codimension two subspaces of $\Zp^n$.  For $n \geq m+1$, Lemma \ref{htrelations} gives $H_p^c(n,m) \supseteq H_p^c(m+1,m) = H_p(m+1,1)$.  So it suffices to analyze $\frac{1}{p} H_p(m+1,1)$, the range of the Nathanson height function on lines in $\mathbb{Z}_p^{m+1}$.
	
	Since we seek to understand values of the form $\frac{1}{p}h_p(V)$, which we term the \emph{fractional height} of $V$, it is helpful to consider the standard embedding of $\mathbb{Z}_p^n$ in the $n$-torus $\mathbb{R}^n/\mathbb{Z}^n$, which identifies $(\bar{a}_1,\dots,\bar{a}_n) \in \mathbb{Z}_p^n$ with $(\frac{a_1}{p}, \dots, \frac{a_n}{p}) \in \mathbb{R}^n/\mathbb{Z}^n$.  If $V$ is the line generated by $(\bar{a}_1,\dots,\bar{a}_n) \in \mathbb{Z}_p^n$, we write $h_p(\bar{a}_1,\dots,\bar{a}_n)$ for $h_p(V)$.  Then the fractional height of $V$ is
\begin{align*}
	\frac{h_p(\bar{a}_1,\dots,\bar{a}_n)}{p} &= \min \left \{ \frac{h_p(v)}{p} \, \Big \vert \, v = (k\bar{a}_1,\dots,k\bar{a}_n), 1 \leq k \leq p-1 \right \} \\
	&= \min \left \{ \frac{a_1^\prime}{p} + \cdots + \frac{a_n^\prime}{p} \, \Big \vert \, (\bar{a}_1^\prime,\dots,\bar{a}_n^\prime) = (k\bar{a}_1,\dots,k\bar{a}_n), 0 \leq a_i^\prime < p, 1 \leq k \leq p-1 \right \} \\
		    &= \min \bigg \{ \left \{\frac{k a_1}{p} \right \} + \cdots + \left \{\frac{k a_n}{p} \right \} \, \Big \vert \, 1 \leq k \leq p-1\bigg\},
\end{align*}
where $\{y\}$ denotes the fractional part of $y \in \mathbb{R}$.  If $x \in \mathbb{R}$ and $(b_1,\dots,b_n) \in \mathbb{Z}^n$, we write
$$
f_x(b_1,\dots,b_n) =  \liminf_{t\rightarrow x} \, \{b_1 t\} + \dots + \{b_n t\}.
$$
Note that when $b_i x \notin \mathbb{Z}$ for every $i$, we have $f_x(b_1,\dots,b_n) = \{b_1 x\} + \dots + \{b_n x\}.$  Then for any $(a_1,\dots,a_n) \in \mathbb{Z}^n$ with $-p < a_i < p$, we have
\begin{equation}\label{feqn}
\frac{h_p(\bar{a}_1,\dots,\bar{a}_n)}{p}= \min_{1\leq k \leq p-1} f_{k/p}(a_1,\dots,a_n).
\end{equation}

The $\liminf$ in the definition of $f_x$ rules out artificially low values when coordinates of opposite signs have fractional part $0$ simultaneously.  For example, $f_{1/2}(6,-2,-3) = 3/2$, not $\{3\}+\{-1\}+\{-3/2\} = 1/2$.  If $p$ grows large and $(a_1,\dots,a_n)$ is fixed, $\{ \frac{1}{p},\dots,\frac{p-1}{p} \}$ becomes dense in $[0,1)$ and $\frac{1}{p}h_p(\bar{a}_1,\dots,\bar{a}_n)$ approaches the global minimum of $f_x(a_1,\dots,a_n)$.  This suggests the following definition.

\begin{definition}
Let $a_1,\dots,a_n \in \mathbb{Z}$.  Then define
	$$h_\infty(a_1,\dots,a_n) =  \min_{x \in [0,1)} f_x(a_1,\dots,a_n).$$
\end{definition}
For any prime $p$ and any $(a_1,\dots,a_n) \in \mathbb{Z}^n$ with $-p < a_i < p$, we have
$$h_\infty(a_1,\dots,a_n) =  \min_{x \in [0,1]} f_x(a_1,\dots,a_n) \leq  \min_{1\leq k \leq p-1} f_{k/p}(a_1,\dots,a_n) = \frac{1}{p} h_p(\bar{a}_1,\dots,\bar{a}_n).$$

Each integral representative of a point in $V \subset \mathbb{Z}_p^n$ whose coordinates have absolute value less than $p$ gives a lower bound on $h_p(V)$.  If $a_1,\dots, a_n \geq 0$, then $h_\infty(a_1,\dots,a_n) = 0$, but if some of the $a_i$'s have opposite signs, then $h_\infty > 0$ and the bound is nontrivial.  This raises two natural questions: how can we calculate the lower bound, and how tight is it?

The following trivial properties of $f_x$ help us to answer those questions.  
\begin{remark}	Let $a_1,\dots,a_n \in \mathbb{Z}$.
\begin{enumerate}[(a)]
\item The function $f_x(a_1,\dots,a_n)$ is piecewise linear in $x$ with slope $a_1+\dots+a_n$.
\item The discontinuities of $f_x(a_1,\dots,a_n)$ occur at rational numbers with denominators $|a_i|$.
\item The local minima of $f_x(a_1,\dots,a_n)$ occur at rational numbers with denominators $|a_i|$.
\end{enumerate}
\end{remark}

While $h_\infty(a_1,\dots,a_n)$ can be difficult to compute in general, we can handle some simple cases easily.

\begin{lemma}\label{somehinfs}
\begin{enumerate}[(a)]
\item Let $a$ and $b$ be integers with $\gcd(a,b) = g$.  Then $h_\infty(a,b) = 0$ if $ab \geq 0$ and $h_\infty(a,b) = \min(\frac{g}{|a|},\frac{g}{|b|})$ otherwise.
\item Let $p$ be a prime and let $a$ and $b$ be positive integers such that $a+b \leq p$.  Then $h_\infty(-p,a,b) = \frac{1}{p} h_p(\bar{a},\bar{b})$.
\end{enumerate}
\end{lemma}
\begin{proof}
\begin{enumerate}[(a)]
\item By property (c) above, $f_x(a,b)$ obtains its minimum value when $x = 0$ or when $x$ is a rational number with denominator $|a|$ or $|b|$.  If $ab \geq 0$, then $f_0(a,b) = 0$, so $h_\infty(a,b) = 0$.  If $ab < 0$, we first consider the case $x = \frac{k}{|a|}$ for some $k > 0$.  We have 
$$f_\frac{k}{|a|}(a,b) = \left \{\frac{ka}{|a|} \right \} + \left \{\frac{kb}{|a|} \right \} =  \left \{\frac{kb}{|a|} \right \} \geq \frac{g}{|a|},$$
with equality when $k$ is the multiplicative inverse of $\frac{b}{g}$ modulo $|a|$.  Similarly, $f_\frac{k}{|b|}(a,b) \geq \frac{g}{|b|}$, with equality $k$ is the multiplicative inverse of $\frac{a}{g}$ modulo $|b|$.  Finally, $f_0(a,b) = 1$.  Thus $h_\infty(a,b) = \min(\frac{g}{|a|},\frac{g}{|b|})$.

\item If $a+b = p$, then $f_x(-p,a,b) = 1$ for all $x$.  If $a+b< p$, then $f_x(-p,a,b)$ is an decreasing function of $x$ where it is continuous.  In fact, it is decreasing as long as $\{-px\} \neq 0$, since $f_x(-p,a,b)$ decreases by $1$ when $x$ crosses $\frac{k}{a}$ or $\frac{k}{b}$.  Thus the minimum of $f_x(-p,a,b)$ occurs when $x = \frac{k}{p}$ for some $k\geq 0$.  Note that
$$\min_{1\leq k \leq p-1} f_\frac{k}{p}(-p,a,b) = \min_{1\leq k \leq p-1} \left \{\frac{ka}{p} \right \} + \left \{\frac{kb}{p} \right \} = \frac{1}{p} h_p(\bar{a},\bar{b}).$$
Since $f_0(-p,a,b) = 1 \geq \frac{1}{p} h_p(\bar{a},\bar{b})$, we have $h_\infty(-p,a,b) = \frac{1}{p} h_p(\bar{a},\bar{b})$. \qedhere 
\end{enumerate} \end{proof}

We can use Lemma \ref{somehinfs} to describe the range of $h_\infty$ when $n = 1,2,$ or $3$.

\begin{lemma}\label{hinfrange}
For any integer $n \geq 2$, let $H_\infty(n) = \{h_\infty(a_1,\dots,a_n) \mid (a_1,\dots,a_n) \in \mathbb{Z}^n \setminus \{0\} \}$.  Then
\begin{enumerate}[(a)]
\item $H_\infty(1) = \{0\};$
\item $H_\infty(2) = \{0\} \cup \{\frac{1}{b}\mid b \in \mathbb{Z}^+\};$
\item $H_\infty(3) \supset \frac{1}{p} H_p(2,1)$ for every prime $p$.
\end{enumerate}
\end{lemma}
\begin{proof}
The first statement is trivial, and the second follows trivially from Lemma \ref{somehinfs}(a).  Now let $V \in \Grp(2,1)$ for some prime $p$.  Since $h_p(V) \leq p$, there are positive integers $a$ and $b$ with $a+b \leq p$ and $(\bar{a},\bar{b}) \in V$.  Lemma \ref{somehinfs} gives $h_\infty(-p,a,b) =  \frac{1}{p}h_p(\bar{a},\bar{b}) = \frac{1}{p}h_p(V)$.
\end{proof}

We can also use the above properties of $f_x$ to show that $h_\infty$ is a tight lower bound on $h_p$ when the $a_i$ are small in proportion to $p$.

\begin{lemma}\label{hinflemma}
Let n be a positive integer, $(a_1,\dots,a_n) \in \mathbb{Z}^n \setminus \{0\}$, and $p$ be a prime such that $p > \max_{i \neq j} |a_ia_j| - \min_i |a_i|$ and $p > \max_i |a_i|$.  Then 
$$0 \leq \frac{1}{p} h_p(\bar{a}_1,\dots,\bar{a}_n) - h_\infty(a_1,\dots,a_n) \leq \frac{|a_1+\dots+a_n|}{p}.$$
\end{lemma}

\begin{proof}
We first consider the case where all but one of the $a_i$ is $0$.  Then $h_\infty(a_1,\dots,a_n) = 0$ and $h_p(\bar{a}_1,\dots,\bar{a}_n) = 1 \leq |a_1+\dots+a_n|$, so the inequality holds.  

Now suppose at least two of the $a_i$ are nonzero.  Let $f_x(a_1,\dots,a_n)$ obtain its minimum value at $x_0 = \frac{k}{|a_i|}$.  By definition, $f_x(a_1,\dots,a_n)$ is either left or right continuous at $x_0$.  Suppose it is left continuous (the proof is similar in the other case).  Let $x_1 = \frac{k^\prime}{|a_j|}$ be the least rational number greater than $x_0$ with denominator in $\{|a_1|,\dots,|a_n|\}$, so $f_x(a_1,\dots,a_n)$ is continuous on $[x_0,x_1)$.  Let $x_p$ be the least rational number greater than $x_0$ with denominator $p$.  We consider two cases.  If $i = j$, then
$$x_p \in \left [x_0,x_0+\frac{1}{p} \right] \subset \left [x_0,x_0+\frac{1}{|a_i|} \right ) =  \left [x_0,x_1 \right ),$$
since $p > |a_i|$.

If $i \neq j$, then
  $$x_p \in \left [x_0,x_0+\frac{1}{p} \right] \subset \left [x_0,x_0+\frac{1}{|a_ia_j|}+\frac{1}{p|a_j|} \right )
	\subseteq  \left [x_0,x_1+\frac{1}{p|a_j|} \right ),$$
since $p > |a_ia_j| - |a_i|$ and $x_1 - x_0 \geq \frac{1}{|a_ia_j|}$.  But $x_p$ and $x_1$ have denominators $p$ and $|a_j|$, respectively, so $x_p \notin \left (x_1 - \frac{1}{p|a_j|}, x_1 + \frac{1}{p|a_j|} \right )$.  Thus, $x_p \in \left [x_0, x_1 - \frac{1}{p|a_j|} \right ] \subset [x_0,x_1)$.

By properties (a) and (b) of $f_x$ above, $f_x(a_1,\dots,a_n)$ is linear on $[x_0,x_1)$ with slope $a_1 + \dots + a_n$, so 
$$f_{x_p}(a_1,\dots,a_n) - h_\infty(a_1,\dots,a_n) \leq |a_1+\dots+a_n|(x_p-x_0) \leq \frac{|a_1+\dots+a_n|}{p}.$$
But 
$$h_\infty(a_1,\dots,a_n)  \leq \frac{1}{p} h_p(\bar{a}_1,\dots,\bar{a}_n) \leq f_{x_p}(a_1,\dots,a_n),$$
which yields the desired inequality.
\end{proof}

We can use Lemmas \ref{somehinfs} and \ref{hinflemma} to give an alternate proof of Corollary \ref{nathansoncor}, and hence an alternate proof of Nathanson's conjecture.

\begin{proof}[Alternate proof of Corollary \ref{nathansoncor}]
	Let $V$ be a line in $\mathbb{Z}_p^2$.  By Lemma \ref{pigeonlemma}, we can choose $a_1,a_2 \in [-\lceil \sqrt{p} \rceil, \lceil \sqrt{p} \rceil ]$ such that $V$ is generated by $(\bar{a}_1,\bar{a}_2) \in \mathbb{Z}_p^2$.  Since multiplication by the inverse of $\gcd(a_1,a_2)$ modulo $p$ fixes $V$, we may assume that $\gcd(a_1,a_2) = 1$.  Then $p > |a_1a_2| - \min(|a_1|,|a_2|)$ and $p > \max(|a_1|,|a_2|)$,
so by Lemma \ref{hinflemma},
	$$0 \leq \frac{1}{p} h_p(\bar{a}_1,\bar{a}_2) - h_\infty(a_1,a_2) < \frac{|a_1+a_2|}{p}.$$
If $a_1a_2 \geq 0$, Lemma \ref{somehinfs} gives $h_\infty(a_1,a_2) = 0$.  In that case,
	$$0 < \frac{1}{p} h_p(V) \leq \frac{|a_1+a_2|}{p} < \frac{2}{\sqrt{p}} + \frac{1}{p},$$
so $0 \leq h_p(V) < 2\sqrt{p} + 1.$

If $a_1a_2<0$, Lemma \ref{somehinfs} gives $h_\infty(a_1,a_2) = \min(\frac{1}{|a_1|}, \frac{1}{|a_2|})$.  Set $b = \max(|a_1|,|a_2|) \leq \left \lceil \sqrt{p} \right \rceil$ and note $|a_1+a_2| < \max(|a_1|,|a_2|) = b$.  Then we have
$$0 \leq \frac{h_p(V)}{p} - \frac{1}{b} < \frac{b}{p},$$
so $0 \leq h_p(V) - p/b < b$.
 \end{proof}

When $V$ is a line in $\mathbb{Z}_p^n$ for $n>2$, there are not necessarily any nonzero points in $V$ whose coordinates satisfy $\max_{i \neq j} |a_ia_j| - \min_i |a_i| < p$, so Lemma \ref{hinflemma} does not supply a complete description of $H_p(n,1)$.  However, if we fix $(a_1,\dots,a_n) \in \mathbb{Z}^n \setminus \{0\}$, then for any prime $p > \max_{i \neq j} |a_ia_j|$, we have 
$$0 \leq \frac{1}{p} h_p(a_1,\dots,a_n) - h_\infty(a_1,\dots,a_n) < \frac{n}{\sqrt{p}}.$$  
Then (taking $n = m+1$) $H_\infty(m+1) \subseteq \limsup_{q \rightarrow \infty} \frac{1}{q} H_q(m+1,1)$.  Lemma \ref{htrelations} states that $H_q(m+1,1) = H^c_q(m+1,m) \subseteq H^c_q(n,m)$ for any $n\geq m+1$, so we have the following theorem.

\begin{theorem}
Let $1 \leq m < n$ be integers.  Then 
$$H_\infty(m+1) \subseteq \limsup_{q \rightarrow \infty} \frac{1}{q} H_q^c(n,m)$$
\end{theorem}

Combining the above theorem with Lemma \ref{hinfrange} gives:
\begin{theorem}  Let $n \geq 3$.  Then
$$\bigcup_{p} \frac{1}{p} H_p(2,1) \subseteq \limsup_{q \rightarrow \infty} \frac{1}{q} H_q^c(2,n). $$
\end{theorem}

%%%%%
\section{Almost all subspaces have small Nathanson height}

In the previous section, we showed that there exist subspaces of $\mathbb{Z}_p^n$ with Nathanson heights that are large in proportion to $p$.  Here, we show that such subspaces are rare.

\begin{theorem}\label{smallheightsn}
Let $\varepsilon > 0$, let $n \geq 2$, and let  $1\leq m \leq n$.  The fraction of $m$-dimensional subspaces of $\mathbb{Z}_p^n$ with Nathanson height greater than $\varepsilon p$ is less than $c_{n,\varepsilon}p^{-1/(n-1)}$ where $c_{n,\varepsilon}$ is a positive constant depending only on $n$ and $\varepsilon$.  In particular, as $p \rightarrow \infty$, this fraction goes to zero.
\end{theorem}

The proof is geometric: for an arbitrary subspace $V$ of $\mathbb{Z}_p^n$, we consider the lattice of all representatives in $\mathbb{Z}^n$ of elements of $V$, and look at a box around the origin large enough that it contains a sizable set of lattice points.  If this set has full rank, then at least one of the points lies in the positive quadrant, so $h_p(V)$ is small.  Otherwise, we show that $V$ must intersect some smaller box, but the fraction of such subspaces is less than $c_{n,\varepsilon}p^{-1/(n-1)}$.  We begin with a lemma that is useful in handling the full-rank case.

\begin{lemma}\label{rickylemma}
Let $v_1,\dots,v_n$ be a basis of a lattice in $\mathbb{R}^n$ with $|v_i| \leq 1$.  Then any ball in $\mathbb{R}^n$ with radius $r$ at least $\frac{\sqrt{n}}{2}$ contains a point of the lattice.
\end{lemma}

\begin{proof}
We proceed by induction.  In the case $n=1$, the lattice is just an arithmetic progression with common difference at most 1.  A ball of radius $1/2$ is a line segment of length 1, so it clearly contains one of the lattice points.

Now suppose that the statement holds for $n$, and let $v_1,\dots,v_{n+1}$ be a basis of a lattice $L$ in $\mathbb{R}^{n+1}$ with $|v_i| \leq 1$.  Let $S$ be a ball in $\mathbb{R}^{n+1}$ of radius at least $ \frac{\sqrt{n+1}}{2}$.  Consider the set of $n$-dimensional cosets of a sublattice of the form $L_m = \{m v_{n+1} + a_1 v_1 + \dots + a_n v_n\mid a_i \in \mathbb{Z} \}$, for each integer $m$.  Each $L_m$ is contained in an $n$-dimensional hyperplane $H_m$, and $H_m$ and $H_{m+1}$ are separated by perpendicular distance at most $|v_{n+1}| \leq 1$.  Therefore, the distance from the center of $S$ to some $H_{m_0}$ is at most $1/2$.  Then $S^\prime = S \cap H_{m_0}$ is a ball in $H_{m_0}$ with radius at least $$\sqrt{\left (\frac{\sqrt{n+1}}{2} \right)^2 - \left (\frac{1}{2} \right )^2} = \frac{\sqrt{n}}{2}.$$  If we consider the point $m v_{n+1}$ as the origin of $\mathbb{R}^n$ embedded in $\mathbb{R}^{n+1}$, the lattice $L_m$ is generated by $v_1,\dots,v_n$.  By the inductive hypothesis, $S^\prime$ contains a point $v$ of $L_m$.  Then $v \in L_m \subset L$ and $S^\prime \subset S$ implies that $v \in S\cap L$, completing the induction.
\end{proof}

The next lemma is our main tool for proving Theorem \ref{smallheightsn}.  We show that any line $V \in \mathbb{Z}_p^n$ either intersects a box with side length on the order of $t p^{ (n-1)/n}$ in the positive quadrant of $\mathbb{Z}^n$, and hence has small Nathanson height, or intersects a box with side length approximately $t^{-1/(n-1)}p^{(n-1)/n}$ centered at the origin.  The negative exponent of $t$ in the latter expression allows us to shrink the size of the box around the origin as $p \rightarrow \infty$ while maintaining the side length of the first box as a fraction of $p$.

\begin{lemma}\label{bigboxlemma}
 Let $n \geq 2$ and $t \geq 2$ be integers.  Let $V$ be a line in $\Zp^n$.  Then there exists some nonzero $(\bar{a}_1,\dots,\bar{a}_n) \in V$ such that  either $0 \leq a_i \leq 2 n t p^{ (n-1)/n}  \mbox{ for } 1 \leq i \leq n$ or $0 \leq |\bar{a}_i|_p \leq 8nt^{-1/(n-1)}p^{(n-1)/n} \mbox{ for } 1 \leq i \leq n$.
\end{lemma}
\begin{proof}
Let $s = \left \lceil t p^{ (n-1)/n} \right \rceil $.  Let $q = \left \lfloor \frac{p}{s} \right \rfloor$.  We can then partition $\mathbb{Z}_p$ into $q$ intervals of length $s$ or $s+1$.  This division of $\mathbb{Z}_p$ into $q$ intervals yields a division of $\mathbb{Z}_p^n$ into $q^{n}$ $n$-dimensional boxes with sides $s$ or $s+1$.  We have 
$$q^{n} = \left \lfloor \frac{p}{s} \right \rfloor^{n} < \left (\frac{p}{s} \right )^{n} < \left (\frac{p}{tp^{ (n-1)/n}}\right )^{n} = \frac{p}{t^n}.$$
	Since $V$ contains $p$ points, some box contains at least 
$$\left \lceil \frac{p}{q^n}\right \rceil \geq \frac{p}{p/t^n} + 1 = t^n + 1$$
points of $V$.  Label them $v_0,v_1,\dots,v_{t^n}$.  Then for $0\leq i \leq t^n$, the points $v_i^\prime = v_i - v_0$ are distinct elements of $V$ lying in $B_s = [-s,s]^n \subset \mathbb{R}^n$.
	
	Now consider the $v_i^\prime$ as points in $\mathbb{Z}^n$ with coordinates whose absolute values are at most $s$.  Suppose that they do not lie in any $(n-1)$-dimensional subspace.  Then some subset of size $n$, say $v_1^\prime,\dots, v_{n}^\prime$, forms the basis of an $n$-dimensional lattice in $\mathbb{Z}^n$.  Set $w_i = \frac{1}{s \sqrt{n}}v_i \in \mathbb{R}^{n}$.  Then $w_1,\dots,w_n$ form the basis of a lattice in $\mathbb{R}^{n}$, and $||w_i||^2 = \frac{||v_i||^2}{s^2n} < \frac{s^2n}{s^2n} = 1$.  Let $S$ be the ball with center $(\frac{\sqrt{n}}{2},\dots,\frac{\sqrt{n}}{2})$ and radius $\frac{\sqrt{n}}{2}$.  By Lemma \ref{rickylemma}, $S$ contains some point on the lattice generated by the $w_i$, say $w = c_1 w_1 + \cdots + c_{n} w_{n}$.  We have 
	$$(a_1,\dots,a_n) = c_1 v_1 + \cdots +  c_nv_n	= (s \sqrt{n}) w \in s \sqrt{n} \cdot (S \cap V).$$
	Furthermore, $s \sqrt{n}S$ has center $(\frac{sn}{2},\dots,\frac{sn}{2})$ and radius $\frac{sn}{2}$, so it is contained in $[0,sn]^n \subset \mathbb{R}^n$.  Thus $a_1,\dots,a_n \in \{0,1,\dots,sn\}$, so $0 \leq a_i \leq 2 n t p^{ (n-1)/n}$ as desired.

	On the other hand, suppose that the lattice $L$ generated by $v_0^\prime,v_1^\prime,\dots,v_{t^n}^\prime$ has rank $m<n$.  Then $L$ is contained in some $m$-dimensional subspace $D$ of $\mathbb{R}^n$.  Let $B_l = [-l,l]^n \subset \mathbb{R}^n$ be the $n$-dimensional box of side length $2l$ centered at the origin, let $B_s^{\prime} = B_s \cap D$, and let $B_l^{\prime} = B_l \cap D$---note that both $B_s^{\prime}$ and $B_l^{\prime}$ contain the origin and are symmetric about it.  We wish to find $l$ large enough that $B_l^\prime$ contains some non-zero lattice point.  By Minkowski's Theorem on the Geometry of Numbers, it suffices to show that $\mbox{vol}(B_l^\prime) \geq 2^m \mbox{vol}(L)$ (see \cite{T1}, p. 134, for a proof).
	
	Note $L\cap B_s^\prime$ contains at least $t^n+1$ points, since $0, v_1^\prime,\dots,v_{t^n}^\prime \in L\cap B_s^\prime$.  We further show that $\#(L\cap B_s^\prime)$ is bounded above by $2(m+1)! \cdot \frac{\vol(B_s^\prime)}{\vol(L)}$: pick a linear function $f: \mathbb{R}^n \rightarrow \mathbb{R}$ such that $f \vert_L$ is injective, and order the points of $L\cap B_s^\prime$ according to their image under $f$.  Take each block of $m+1$ consecutive points and form a simplex.  The resulting $\left \lfloor \frac{\#(L\cap B_s^\prime)}{m+1} \right \rfloor > \frac{\#(L\cap B_s^\prime)}{2(m+1)}$ simplices are disjoint and contained in the convex hull of $L\cap B_s^\prime$.  Finally, each simplex has volume at least $\frac{\vol(L)}{m!}$.  Combining these inequalities gives
$$ \frac{\vol(B_l^\prime)}{\vol(L)} = \frac{\vol(B_l^\prime)}{\vol(B_s^\prime)} \cdot \frac{\vol(B_s^\prime)}{\vol(L)}
	\geq \left (\frac{l}{s} \right )^m \frac{\#(L\cap B_s^\prime)}{2(m+1)!}
	\geq \left (\frac{l}{s} \right )^m \frac{t^n+1}{2(m+1)!}. $$
Let $l = 8nt^{-1/(n-1)}p^{(n-1)/n} > 4nst^{-n/(n-1)}$.  Then
$$\frac{\vol(B_l^\prime)}{\vol(L)} \geq \left (\frac{l}{s} \right )^m \frac{t^n+1}{2(m+1)!} 
\geq \left (\frac{4n}{t^\frac{n}{n-1}} \right )^m\frac{t^n+1}{2(m+1)!}
\geq  \left (\frac{4^m(m+1)!}{t^n} \right )\frac{t^n+1}{2(m+1)!}
\geq 2^m,$$
since $n-1 \geq m$.  Therefore, some nonzero point of $L \subset V$ lies in $B_l^\prime \subset B_l$, and our statement follows.
\end{proof}

We can use Lemma \ref{bigboxlemma} to prove Theorem \ref{smallheightsn} in the case where $V$ is a one dimensional subspace of $\mathbb{Z}_p^n$.

\begin{lemma}\label{smallheights}
Fix $\varepsilon > 0$.  The fraction of lines in $\mathbb{Z}_p^n$ with Nathanson height greater than $\varepsilon p$ is less than $c_n \varepsilon^{-n/(n-1)}p^{-1/(n-1)}$ where $c_n$ is a positive constant depending only on $n$.  In particular, as $p \rightarrow \infty$, this fraction goes to zero.
\end{lemma}
\begin{proof}
	Set $t = \lfloor \varepsilon p^{1/n}/4n^2 \rfloor$ in Lemma \ref{bigboxlemma}.  Then any line $V \in \mathbb{Z}_p^n$ contains some point $(\bar{a}_1,\dots,\bar{a}_n)$ with either $0 \leq a_i \leq \varepsilon p/2n < \varepsilon p/n$ or $0 \leq |\bar{a}_i|_p \leq l$, where $l = 8nt^{-1/(n-1)}p^{(n-1)/n}.$  If $V$ satisfies the first condition, then $h_p(V) < \varepsilon p$.  Since there are only $(2l)^n$ points in $\mathbb{Z}_p^n$ satisfying the second condition, there are at most $(2l)^n$ lines in $\mathbb{Z}_p^n$ with Nathanson height greater than $\varepsilon p$, out of $(p^n - 1)/(p-1)$ total lines.  Taking this ratio and reducing gives the desired result.
\end{proof}

Since every subspace of dimension greater than one contains many lines, we can use Lemma \ref{smallheights} to prove the general case of Theorem \ref{smallheightsn}.

\begin{proof}[Proof of Theorem \ref{smallheightsn}]
We proceed with a probabilistic argument.  Choose a random $n$-dimensional subspace  $V$ of $\mathbb{Z}_p^n$ uniformly, and choose a random line $V^\prime \subset V$ uniformly.  Then $V^\prime$ has uniform distribution over the set of lines in $\mathbb{Z}_p^n$, so by Lemma \ref{smallheights}, the probability that it has Nathanson height greater than $\varepsilon p$ is less than $c_{n,\varepsilon}p^{-1/(n-1)}$, where $c_{n,\varepsilon} = c_n \varepsilon^{-n/(n-1)}$.  The Nathanson height of $V$ is bounded above by the Nathanson height of its subspaces, so $V$ has Nathanson height greater than $\varepsilon p$ with lower probability than $V^\prime$.  Therefore the fraction of $m$-dimensional subspaces of $\mathbb{Z}_p^n$ with Nathanson height greater than $\varepsilon p$ is less than $c_{n,\varepsilon}p^{-1/(n-1)}$.
\end{proof}

\section{Conclusion and open problems}

The behavior of the Nathanson height function is well understood for subspaces of codimension one.  As $p$ grows large, the fractional heights approach the set $\Bset = H_\infty(2)$, with $\Dist(H_\infty(2),\frac{1}{p} H_p^c(n,1)) \leq O(p^{-1/n})$.  For higher codimension, we could only show inclusion: 
$$H_\infty(m+1) \subseteq \limsup_{q \rightarrow \infty} \frac{1}{q} H_q^c(n,m).$$  
We conjecture, though, that we have the same equality and asymptotics as in the codimension one case.
\begin{conjecture}
Let $1 \leq m <n$.  Then
$$\Dist(H_\infty(m+1),\frac{1}{p} H_p^c(n,m)) \leq O(p^{-1/n}).$$
\end{conjecture}

The properties of the set $H_\infty(m)$ are also well understood only for $m = 2$.  Finally, the properties of the coheight and width functions are largely unknown for finite abelian groups other than $\Zp$.   We conjecture, though, that analogues of Theorems \ref{coheightthm} and \ref{widththm} hold when $G$ is an arbitrary cyclic group.

\section{Acknowledgments}
This research was done at the University of Minnesota Duluth with the financial support of the National Science Foundation (grant number DMS-00447070-001) and the National Security Agency (grant number H98230-06-1-0013).  I would like to thank Ricky Liu and Reid Barton for their helpful suggestions and careful editing.  I would also like to thank Joe Gallian for suggesting the topic and running the summer research program.

\end{document}